\documentclass[12pt, amscd]{amsart}
\usepackage{amscd}
\usepackage{verbatim}

\input xy
\xyoption{all}

\usepackage{amssymb}

\DeclareMathAlphabet{\cat}{OT1}{cmss}{m}{sl}

\newtheorem{theorem}{Theorem}[section]

\newtheorem{proposition}[theorem]{Proposition}
\newtheorem{lemma}[theorem]{Lemma}
\newtheorem{corollary}[theorem]{Corollary}

\theoremstyle{definition}
\newtheorem{remark}[theorem]{Remark}

\newcommand{\xra}{\xrightarrow}

\newcommand{\tens}{\otimes}
\newcommand{\inv}{^{-1}}
\newcommand{\gmu}{\boldsymbol{\mu}}

\newcommand{\ind}{\operatorname{\hspace{0.3mm}ind}}
\newcommand{\ch}{\operatorname{char}}

\newcommand{\Br}{\operatorname{Br}}
\newcommand{\Spec}{\operatorname{Spec}}

\newcommand{\gGL}{\operatorname{\mathbf{GL}}}
\newcommand{\ga}{\operatorname{\mathbb{G}}_a}

\newcommand{\ed}{\operatorname{ed}}

\newcommand{\td}{\operatorname{tr.deg}}

\newcommand{\A}{\mathbb{A}}

\newcommand{\Z}{\mathbb{Z}}

\newcommand{\cB}{\mathcal B}

\usepackage[hypertex]{hyperref}

\title[Essential dimension of simple algebras in positive characteristic] 
{Essential dimension of simple algebras in positive characteristic}

\author
[S. Baek] {Sanghoon Baek}

\address
{Department of Mathematics and Statistics, University of Ottawa, Canada}

\email {sbaek@uottawa.ca}

\thanks{The work has been supported by NSERC
Discovery 385795-2010, Accelerator Supplement 396100-2010 grants, and Erhard Neher's NSERC Discovery grant.}

\begin{document}

\begin{abstract}
Let $p$ be a prime integer, $1\leq s\leq r$ integers, $F$ a field of characteristic $p$. Let $\cat{Dec}_{p^r}$ denote the class of the tensor product of $r$ $p$-symbols and $\cat{Alg}_{p^r,p^s}$ denote the class of central simple algebras of degree $p^r$ and exponent dividing $p^s$. For any integers $s<r$, we find a lower bound for the essential $p$-dimension of $\cat{Alg}_{p^r,p^s}$. Furthermore, we compute upper bounds for $\cat{Dec}_{p^r}$ and $\cat{Alg}_{8,2}$ over $\ch(F)=p$ and $\ch(F)=2$, respectively. As a result, we show $\ed_{2}(\cat{Alg}_{4,2})=\ed(\cat{Alg}_{4,2})=\ed_{2}(\gGL_{4}/\gmu_{2})=\ed(\gGL_{4}/\gmu_{2})=3$ and $3\leq \ed(\cat{Alg}_{8,2})=\ed(\gGL_{8}/\gmu_{2})\leq 10$ over a field of characteristic $2$.
\end{abstract}

\maketitle

\section{Introduction}

A numerical invariant, essential dimension of algebraic groups was introduced by Reichstein and was generalized to algebraic structures by Merkurjev. We refer to \cite{Merkurjev09} for the definition of essential dimension and denote by $\ed$ and $\ed_{p}$ the essential dimension and essential $p$-dimension, respectively.

Let $F$ be a field, $\cat{Fields}/F$ the category of field extensions over $F$, and $\cat{Sets}$ the category of sets. For every integer $n\geq 1$, a divisor $m$ of $n$ and any field extension $E/F$, let
\[\cat{Alg}_{n,m}: \cat{Fields}/F\to\cat{Sets}\] be the functor taking a field extension $K/F$ to the set of isomorphism classes of central simple $K$-algebras of degree $n$ and exponent dividing $m$. Then, there is a natural bijection between $H^1(K,\gGL_{n}/\gmu_{m})$ and $\cat{Alg}_{n,m}(K)$ (see \cite[Example 1.1]{BM09}), thus $\ed(\cat{Alg}_{n,m})=\ed(\gGL_{n}/\gmu_{m})$ and $\ed_{p}(\cat{Alg}_{n,m})=\ed_{p}(\gGL_{n}/\gmu_{m})$.

Let $F$ be a field of characteristic $p$. For  $a\in F$ and $b\in F^{\times}$, we write $\big[a,b\big)_{p}$ for the central simple algebra over $F$ generated by $u$ and $v$ satisfying $u^{p}-u=a$, $v^{p}=b$ and $vu=uv+v$ (it is called a symbol $p$-algebra). For a field extension $E/F$, let
\[\cat{Dec}_{p^r}: \cat{Fields}/F\to\cat{Sets}\] be the functor taking a field extension $K/F$ to the set of isomorphism classes of the tensor product of $r$ $p$-symbols over $E$.

Some computations of the essential dimension and essential $p$-dimension of $\cat{Alg}_{m,n}$ have been done. But most of them have the restriction $\ch(F)\neq p$ on the base field $F$. In this paper, for any integers $r>s$, we find a new lower bound for $\ed_{p}(\cat{Alg}_{p^r,p^s})$ over $\ch(F)=p$. Moreover, we compute upper bounds for $\cat{Dec}_{p^r}$ and $\cat{Alg}_{8,2}$ over $\ch(F)=p$ and $\ch(F)=2$, respectively. As a result, we get:

\begin{theorem}
Let $F$ be a field of characteristic $2$. Then \[\ed_{2}(\cat{Alg}_{4,2})=\ed(\cat{Alg}_{4,2})=\ed_{2}(\gGL_{4}/\gmu_{2})=\ed(\gGL_{4}/\gmu_{2})=3.\]
\end{theorem}
\begin{proof}
The lower bound $3\leq \ed_{2}(\cat{Alg}_{4,2})$ follows from Corollary \ref{maintheorem}. By a Theorem of Albert, we have $\cat{Dec}_{4}=\cat{Alg_{4,2}}$ for $p=2$, thus we get $\ed(\cat{Alg}_{4,2})\leq 3$ by Proposition \ref{decupper}. As $\ed_{2}(\cat{Alg}_{4,2})\leq \ed(\cat{Alg}_{4,2})$, the result follows.
\end{proof}

Corollary \ref{maintheorem} and Corollary \ref{last} give the following:

\begin{theorem}
Let $F$ be a field of characteristic $2$. Then \[3\leq \ed(\cat{Alg}_{8,2})=\ed(\gGL_{8}/\gmu_{2})\leq 10.\]
\end{theorem}

\section{Lower bounds}

\begin{theorem}\label{Tsen}
$($Tsen$)$ Let $K$ be a field of transcendental degree $1$ over an algebraically closed field $F$. Then, for any central division algebra $A$ over $K$, $\ind(A)=\exp(A)=1$, i.e., $A=K$.
\end{theorem}

As an application of Theorem \ref{Tsen}, Reichstein obtained the following result:
\begin{corollary}\label{Reichstein}\cite[Lemma 9.4(a)]{Re00}
Let $F$ be an arbitrary field and $A$ be a division algebra of degree $n\geq 2$. Then $\ed(A)\geq 2$. In particular, for any integers $r$, $s$ and any prime $p$, $\ed_{p}(\cat{Alg}_{p^r,p^s})\geq 2$.
\end{corollary}

Initially, the following theorem is proved under the additional condition that $\ch(F)$ does not divide $\exp(A)$ in \cite{deJong}. In subsequent papers \cite[Theorem 1.0.2]{deJongStarr} and \cite[Theorem 4.2.2.3]{Lieblich}, this condition is removed:
\begin{theorem}\label{periodindex}
$($de Jong$)$ Let $K$ be a field of transcendental degree $2$ over an algebraically closed field $F$. Then, for any central simple algebra $A$ over $K$, $\ind(A)=\exp(A)$.
\end{theorem}

As an application of Theorem \ref{periodindex}, we have the following result:
\begin{corollary}\label{maintheorem}
Let $F$ be an arbitrary field and $p$ be a prime. For any integers $r$ and $s$ with $s<r$, $\ed_{p}(\cat{Alg}_{p^r,p^s})\geq 3$.
\end{corollary}
\begin{proof}
By \cite[Proposition 1.5]{Merkurjev09}, we may replace the base field $F$ by an algebraically closure of $F$. Let $K$ be a field extension of $F$ and $A$ be a central simple algebra over $K$ of $\ind(A)=p^r$ and $\exp(A)|p^s$. Let $E$ be a field extension of $K$ of degree prime to $p$. As $\ind(A)$ is relatively prime to $[E:K]$, we have $\ind(A_{E})=\ind(A)=p^r$. Suppose that $A_{E}\simeq B\tens E$ for some $B\in \cat{Alg}_{p^r,p^s}(L)$ and $\td_{F}(L)=2$. Then, by Theorem \ref{periodindex}, we have $\ind(B)=\exp(B)$. As $p^r=\ind(A_{E})|\ind(B)=\exp(B)$, we get $p^r|\exp(B)$. But it contradicts to $\exp(B)|p^s$. By Corollary \ref{Reichstein}, the result follows.
\end{proof}
\begin{remark}
As we see in \cite[Theorem]{BM10}, the above lower bound $3$ is much less than the best known lower bounds, but these lower bounds are valid only for $\ch(F)\neq p$. Hence, our main application of Corollary \ref{maintheorem} is for the case of $\ch(F)=p$.
\end{remark}

\section{Upper bounds}

\subsection{An upper bound for $\ed(\cat{Dec}_{p^r})$}\label{decupperbd}

\begin{lemma}\label{edzmodpz}\textnormal{\cite[Example 2.3 and page 298]{BerhuyFavi03}}
Let $F$ be a field of characteristic $p$ and $r\geq 1$ be an integer. If $|F|\geq p^{r}$, then $\ed((\Z/p\Z)^{r})=1$.
\end{lemma}
\begin{proof}
From the exact exact sequence \[0\to \Z/p\Z\to \ga\xra{\wp} \ga\to 0,\] we have $H^{1}(E,\Z/p\Z)=E/\wp(E)$ for any field extension $E/F$ where $\wp(x)=x^{p}-x$ for $x\in E$, hence $\ed(\Z/p\Z)=1$. By \cite[Proposition 4.11]{BerhuyFavi03}, we have $\ed((\Z/p\Z)^{r})\geq 1$.

As $|F|\geq p^{r}$, we have an exact sequence \[0\to (\Z/p\Z)^{r}\to \ga\xra{\wp} \ga\to 0.\] It follows from $H^{1}(E,\ga)=0$ that $\ga(E)\to H^{1}(E,(\Z/p\Z)^{r})$ is surjective. By \cite[Proposition 1.3]{Merkurjev09}, we have $\ed((\Z/p\Z)^{r})\leq \dim(\ga)=1$.
\end{proof}

\begin{proposition}\label{decupper}
Let $p$ be a prime integer, $F$ be a field of characteristic $p$. If $|F|\geq p^{r}$, then $\ed(\cat{Dec}_{p^r})\leq r+1$.
\end{proposition}
\begin{proof}
Let \[A=\otimes_{i=1}^{r} \big[a_{i},b_{i}\big)_{p} \in \cat{Dec}_{p^r}(E)\] for a field extension $E/F$. As $\ed((\Z/p\Z)^{r})=1$ by Lemma \ref{edzmodpz}, there exists a sub-extension $E_{0}/F$ of $E/F$ and $c_{i}\in E_{0}$ for all $1\leq i \leq r$ such that $c_{i}\equiv a_{i} \mod \wp(E)$ and $\td_{F}(E_{0})\leq 1$. Therefore, $A$ is defined over $L=E_{0}(b_{1},\cdots,b_{r})$ and $\td_{F}(L)\leq r+1$. Hence, $\ed(A)\leq r+1$ and $\ed(\cat{Dec}_{p^r})\leq r+1$.

\end{proof}

\subsection{An upper bound for $\ed(\cat{Alg}_{8,2})$}\label{charalg2}
In this subsection we assume that the base field $F$ is of characteristic $2$. The upper bound $8$ for $\ed(\cat{Alg}_{8,2})$ over the base field $F$ of characteristic different from $2$ was determined in \cite[Theorem 2.12]{BM09}. We use a similar method to find an upper bound for $\ed(\cat{Alg}_{8,2})$ over the base field $F$ of characteristic $2$.

For a commutative $F$-algebra $R$, $a\in R$ and $b\in R^\times$ we write $\big[a,b\big)_R$ for the quaternion algebra $R\oplus
Ru\oplus Rv\oplus Rw$ with the multiplication table $u^2+u=a, v^2=b,
uv=w=vu+v$. The class of $\big[a,b\big)_R$ in the Brauer group
$\Br(R)$ will be denoted by $\big\{a,b\big\}=\big\{a,b\big\}_R$.

Let $a\in R$ and $S=R[\alpha]:=R[t]/(t^2+t+a)$ with $\alpha^2=\alpha+a$ the quadratic
extension of $R$. We write $N_R(a)$ for the subgroup of $R^\times$
of all element of the form $x^2+xy+ay^2$ with $x,y\in R$, i.e.,
$N_R(a)$ is the image of the norm homomorphism $N_{S/R}:S^\times\to
R^\times$. If $b\in N_R(a)$, then the quaternion algebra
$\big[a,b\big)_R$ is isomorphic to the matrix algebra $M_2(R)$.

\subsubsection{Rowen's construction}

Rowen extended the Tignol's theorem \cite{Tignol78} to a field of characteristic $2$. We recall Rowen's argument in \cite{Rowen}. Let $A$ be a central simple $F$-algebra in $\cat{Alg}_{8,2}(F)$. By \cite{Rowen},
there is a triquadratic splitting extension
$F(\alpha, \beta, \gamma)/F$ of $A$ such that $\alpha^2+\alpha=a, \beta^2+\beta=b,$ and $\gamma^2+ \gamma=c$ for some $a,b,c\in F$.
Let $L=F(\alpha)$. By \cite[Corollary 5]{Rowen}, we have
\begin{equation}\label{overlll}
\big\{A\big\}_L=\big\{b,s\big\}+\big\{c,t\big\}
\end{equation}
in $\Br(L)$ for some $s,t\in L^\times$.

Taking the corestriction for the extension $L/F$ in $(\ref{overlll})$,
we get
\[
0=2\big\{A\big\}=\big\{b,N_{L/F}(s)\big\}+\big\{c,N_{L/F}(t)\big\}
\]
in $\Br(F)$, hence $\big\{b,N_{L/F}(s)\big\}=\big\{c,N_{L/F}(t)\big\}$.
By the chain lemma \cite[Lemma 3]{Rowen}, we have
\[
\big\{b,N_{L/F}(s)\big\}=\big\{d,N_{L/F}(s)\big\}=\big\{d,N_{L/F}(t)\big\}=\big\{c,N_{L/F}(t)\big\}
\]
in $\Br(F)$ for some $d\in F$. Therefore, we get
$\big\{b+d,N_{L/F}(s)\big\}=\big\{c+d,N_{L/F}(t)\big\}=\big\{d,N_{L/F}(st)\big\}=0$. By the proof of \cite[Lemma 2.3]{Dherte95},
\begin{align*}
\big\{b+d,s\big\}&=\big\{b+d,k\big\},\\
\big\{c+d,t\big\}&=\big\{c+d,l\big\},\\
\big\{d,st\big\}&=\big\{d,m\big\}.
\end{align*}
in $\Br(L)$ for some $k,l,m\in F^\times$. It follows from
(\ref{overlll}) that
\[
\big\{A\big\}_L=\big\{b+d,k\big\}_L+\big\{c+d,l\big\}_L+\big\{d,m\big\}_L
\]
in $\Br(L)$. Hence
\begin{equation*}
\big\{A\big\}=\big\{a,e\big\}+\big\{b,k\big\}+\big\{c,l\big\}+\big\{d,klm\big\}
\end{equation*}
in $\Br(F)$ for some $e\in F^\times$.

We shall need the following result:

\begin{lemma}\label{comeeed}
Let $R$ be a commutative $F$-algebra, $a, b\in R$,
$T=R[\alpha]:=R[t]/(t^2+t+a)$ and $x+y\alpha\in T^\times$ such that
$x^2+xy+ay^2=u^2+uv+bv^2$ for some $u,v\in R$. If $v+y\in R^\times$, then
$(v+y)(x+y\alpha)\in N_T(b)$. In particular,
\[
\big\{b,x+y\alpha\big\}_T=\big\{b,v+y\big\}_T.
\]
\end{lemma}
\begin{proof}
The result comes from the following equality
\begin{align*}
(x+y\alpha+u)^2+(x+y\alpha+u)v+bv^2=(x+y\alpha)(v+y).
\end{align*}
\end{proof}

\subsubsection{Classifying Azumaya algebra for $\cat{Alg}_{8,2}$}
Consider the affine space $\A^{13}_F$ with coordinates
$\mathbf{a},\mathbf{b},\mathbf{c},\mathbf{d}, \mathbf{e},\mathbf{u},\mathbf{v},\mathbf{w},\mathbf{x},\mathbf{y},\mathbf{z},\mathbf{m},\mathbf{n}$
and define the rational functions:
\begin{align*}
    \mathbf{f}&=\mathbf{x}\mathbf{z}+\mathbf{w}\mathbf{z}+\mathbf{x}\mathbf{y},\\
    \mathbf{g}&=\mathbf{w}\mathbf{y}+\mathbf{x}\mathbf{z}\mathbf{a},\\
    \mathbf{r}&=(\mathbf{g}^2+\mathbf{g}\mathbf{f}+\mathbf{f}^2\mathbf{a}+\mathbf{m}^2+\mathbf{m}\mathbf{n}),\\
    \mathbf{h}&=(\mathbf{w}^2+\mathbf{w}\mathbf{x}+\mathbf{x}^2\mathbf{a}+1+\mathbf{u}+\mathbf{u}^2\mathbf{d}),\\
    \mathbf{l}&=(\mathbf{y}^2+\mathbf{y}\mathbf{z}+\mathbf{z}^2\mathbf{a}+1+\mathbf{v}+\mathbf{v}^2\mathbf{d}),\\
    \mathbf{p}&=(\mathbf{u}+\mathbf{x})(\mathbf{v}+\mathbf{z})(\mathbf{n}+\mathbf{f}).\\
\end{align*}
Let $X$ be the quasi-affine variety defined by
\begin{equation*}
\mathbf{q}:=\mathbf{a}\mathbf{b}\mathbf{c}\mathbf{d}\mathbf{e}\mathbf{p}(\mathbf{w}^2+\mathbf{w}\mathbf{x}+\mathbf{x}^2\mathbf{a})
(\mathbf{y}^2+\mathbf{y}\mathbf{z}+\mathbf{z}^2\mathbf{a})
(\mathbf{g}^2+\mathbf{g}\mathbf{f}+\mathbf{f}^2\mathbf{a})\neq 0,
\end{equation*}
\begin{equation*}
\mathbf{b}\mathbf{u}^{2}=\mathbf{h},\mathbf{c}\mathbf{v}^{2}=\mathbf{l},\mathbf{d}\mathbf{n}^{2}=\mathbf{r}
\end{equation*}

i.e., $X=\Spec(R)$ with
\[R=F[\mathbf{a},\mathbf{b},\mathbf{c},\mathbf{d},\mathbf{e},\mathbf{u},\mathbf{v},\mathbf{w},\mathbf{x},\mathbf{y},\mathbf{z},\mathbf{m},\mathbf{n},\mathbf{q}\inv]
/<\mathbf{b}\mathbf{u}^{2}+\mathbf{h},\mathbf{c}\mathbf{v}^{2}+\mathbf{l},\mathbf{d}\mathbf{n}^{2}+\mathbf{r}>.\] Let $T=R[\alpha]$ and $S=R[\alpha, \beta, \gamma]$ with $\alpha^2=\alpha+\mathbf{a}, \beta^2=\beta+\mathbf{b}, \gamma^2=\gamma+\mathbf{c}$.
Consider the Azumaya $R$-algebra
\begin{equation}\label{elementtwo}
\cB'=\big[\mathbf{a},\mathbf{e}
\big)_R\tens\big[\mathbf{b},\mathbf{x}+\mathbf{u}\big)_R\tens\big[\mathbf{c},\mathbf{z}+
\mathbf{v}\big)_R\tens\big[\mathbf{d},\mathbf{p}\big)_R.
\end{equation}

By Lemma \ref{comeeed}, we get the followings :
\begin{align*}
(\mathbf{x}+\mathbf{u})(\mathbf{w}+\mathbf{x}\alpha)&\in N_T(\mathbf{b}+\mathbf{d})\subset N_S(\mathbf{d}),\\
(\mathbf{z}+\mathbf{v})(\mathbf{y}+\mathbf{z}\alpha)&\in N_T(\mathbf{c}+\mathbf{d})\subset N_S(\mathbf{d}),\\
(\mathbf{n}+\mathbf{f})(\mathbf{w}+\mathbf{x}\alpha)(\mathbf{y}+\mathbf{z}\alpha)&
\in N_T(\mathbf{d})\subset N_S(\mathbf{d}).\\
\end{align*}
It follows from (\ref{elementtwo}) that
\begin{equation*}
\big\{\cB'\big\}_T=\big\{\mathbf{b},\mathbf{w}+\mathbf{x}\alpha\big\}+\big\{\mathbf{c},\mathbf{y}+\mathbf{z}\alpha\big\}
\end{equation*}
in $\Br(T)$.

Since $\mathbf{p}\in N_S(\mathbf{d})$, $\big[\mathbf{d},\mathbf{p}\big)_S$ is isomorphic to the matrix
algebra $M_2(S)$. In particular,
\[
M_2(R)\subset M_2(S)\simeq
\big[\mathbf{d},\mathbf{p}\big)_S\subset\cB'
\]
and hence $\cB'\simeq M_2(\cB)$ for the centralizer $\cB$ of
$M_2(R)$ in $\cB'$ by the proof of \cite[Theorem 4.4.2]{Herstein}. Then
$\cB$ is an Azumaya $R$-algebra of degree $8$ that is Brauer equivalent
to $\cB'$ by \cite[Theorem 3.10]{Saltman}.

\begin{proposition}\label{classchartwo}
The Azumaya algebra $\cB$ is classifying for $\cat{Alg}_{8,2}$, i.e,
the corresponding $\gGL_8/\gmu_2$-torsor over $X$ is classifying.
\end{proposition}

\begin{proof}
Let $A\in \cat{Alg}_{8,2}(K)$, where $K$ is a field extension of $F$.
We shall find a point $p\in X(K)$ such that $A\simeq \cB(p)$, where $\cB(p):=\cB\tens_{R}K$ with the $F$-algebra homomorphism $R\to K$ given by the point $p$.

Following Rowen's construction, there is a triquadratic splitting
extension $K(\alpha,\beta,\gamma)/K$ of $A$ such that $\alpha^2+\alpha=a, \beta^2+\beta=b,$ and $\gamma^2+ \gamma=c$ for some $a,b,c\in K$. Let $L=K(\alpha)$, so
\[
\big\{A\big\}_{L}=\big\{b,s\big\}+\big\{c,t\big\}
\]
in $\Br(L)$ for some $s=w+x\alpha,$ and $t=y+z\alpha\in
{L}^\times$. We have
\[
\big\{b,w^2+wx+x^2a\big\}=\big\{d,w^2+wx+x^2a\big\}=\big\{d,y^2+yz+z^2a\big\}=\big\{c,y^2+yz+z^2a\big\}
\]
in $\Br(K)$ for some $d\in K$, so $\big\{b+d,w^2+wx+x^2a\big\}=\big\{c+d,y^2+yz+z^2a\big\}=\big\{d,(w^2+wx+x^2a)(y^2+yz+z^2a)\big\}=0$. Hence
\begin{align*}
    w^2+wx+x^2a&={u'}^2+u'u+u^2(b+d),\\
    y^2+yz+z^2a&={v'}^2+v'u+v^2(c+d),\\
    (w^2+wx+x^2a)(y^2+yz+z^2a)&=m^2+mn+n^2d
\end{align*}
for some $u,u',v,v',m,n$ in $K$. Moreover, we may assume that $u'\neq 0$. Replacing $w$, $x$ and $u$ by $wu'$, $xu'$ and $u'u$ respectively, we may assume that $u'=1$. Similarly, we may assume that $v'=1$.

We also may assume that $u\neq x$ by replacing $u$ by $u/(b+d)$. Similarly, we may assume that $v\neq z$ and $n+xz+wz+xy\neq 0$. It follows from Lemma \ref{comeeed} that
\begin{align*}
\big\{b+d,w+x\alpha\big\}&=\big\{b+d,u+x\big\}_{L},\\
\big\{c+d,y+z\alpha\big\}&=\big\{c+d,z+v\big\}_{L},\\
\big\{d,(w+x\alpha)(y+z\alpha)\big\}&=\big\{d,n+xz+wz+xy\big\}_{L}
\end{align*}
in $\Br(L)$. Hence
\[
\big\{A\big\}=\big\{a,e\big\}+\big\{b,u+x\big\}+\big\{c,z+v\big\}+\big\{d,(u+x)(z+v)(n+xz+wz+xy)\big\}
\]
in $\Br(K)$ for some $e\in K^\times$.

Let $p$ be the point $(a,b,c,d,e,u,v,w,x,y,z,m,n)$ in $X(K)$. We have
$\big\{\cB(p)\big\}=\big\{A\big\}$ and hence $\cB(p)\simeq A$ as
$\cB(p)$ and $A$ have the same dimension.
\end{proof}

\begin{corollary}\label{last}
$\ed(\cat{Alg}_{8,2})\leq 10$.
\end{corollary}

\begin{proof}
There is surjective morphism $X\to \cat{Alg}_{8,2}$ by Proposition \ref{classchartwo}. By \cite[Proposition 1.3]{Merkurjev09}, $\ed(\cat{Alg}_{8,2})\leq \dim(X)=10$.
\end{proof}

\end{document}